\documentclass[11pt]{article}
\usepackage{amsmath,amssymb,amsthm,url,palatino}
\usepackage{color}
\usepackage[margin=1.0in,letterpaper]{geometry}
\usepackage{graphicx}
\usepackage{longtable}
\usepackage{tikz}
\usetikzlibrary{shapes,arrows}

\newtheorem{thm}{Theorem}[section]
\newtheorem{lemma}[thm]{Lemma}
\newtheorem{prop}[thm]{Proposition}

\theoremstyle{definition}
\newtheorem{defn}[thm]{Definition}
\newtheorem{example}[thm]{Example}
\newtheorem{algorithm}[thm]{Algorithm}

\DeclareMathOperator{\Sym}{Sym}
\DeclareMathOperator{\Rep}{Rep}
\DeclareMathOperator{\Tw}{Tw}
\DeclareMathOperator{\fix}{fix}
\DeclareMathOperator{\id}{id}
\DeclareMathOperator{\im}{Im}

\renewcommand{\gg}{\mathbf{g}}
\newcommand{\hh}{\mathbf{h}}
\newcommand{\ww}{\mathbf{w}}
\newcommand{\vv}{\mathbf{v}}

\newcommand{\ee}{\mathbf{e}}
\newcommand{\bb}{\mathbf{b}}

\newcommand{\FF}{\mathbb{F}}

\title{Decoding twisted permutation codes}
\author{Robert F.\ Bailey\footnote{School of Science and the Environment (Mathematics), Grenfell Campus, Memorial University, Corner Brook, NL A2H~6P9, Canada. E-mail: \texttt{rbailey@grenfell.mun.ca}} \and Keenan B.\ Nicholson\footnote{Department of Computer Science, Memorial University, St.\ John's, NL A1C~5S7, Canada. E-mail: \texttt{keenanbnicholson@gmail.com}}}

\begin{document}
\maketitle

\begin{abstract}
We consider twisted permutation codes, a class of frequency permutation arrays obtained from finite groups with multiple permutation representations of the same degree, introduced by Gillespie, Praeger and Spiga (and later studied by Akbari, Gillespie and Praeger), and develop a decoding algorithm for such codes based on earlier work of the first author for permutation group codes.  In particular, we show how to implement this algorithm for an infinite family of groups considered by Akbari, Gillespie and Praeger.\\

uncovering-by-bases.\\


\end{abstract}

\section{Introduction} \label{sec:intro}
A {\em permutation code} is an error-correcting code where each codeword is a permutation written in list form (i.e.\ a listing of the symbols from a set of size $n$, where each symbol appears exactly once).  Such a code is also known as a {\em permutation array}, $\mathrm{PA}(n,d)$, where $d$ denotes the minimum Hamming distance.  Permutation codes have a history dating back at least to the 1970s (see~\cite{Blake79}, for instance), but have more recently been considered because of applications including powerline communications~\cite{ChuColbournDukes04}, solid-state memory devices~\cite{Klove10,TamoSchwartz10} and DNA storage of data~\cite{BeeriSchwartz22}.  We note that for two permutations $g,h$ in the symmetric group $S_n$, their Hamming distance is $n-\fix(gh^{-1})$ (where $\fix(g)$ denotes the number of fixed points of $g$).  In the case where the set of permutations forms a subgroup $G$ of $S_n$, the minimum distance is
\[ \min_{ \substack{g\in G \\ g\neq 1} } \left\{ n - \fix(g)\right\}\,\,\,  =  \,\,\, n - \max_{ \substack{g\in G \\ g\neq 1} } \left\{ \fix(g) \right\}. \]
The study of groups of permutations as codes is the subject of several papers of the first author and others~\cite{btubb,ecpg,m12,soc,distenum,thomas}.

A more general notion is that of a {\em frequency permutation array}, $\mathrm{FPA}_\lambda(m,d)$, which is a code of alphabet size $n$ and length $m=\lambda n$, where each codeword contains each of the $n$ symbols exactly $\lambda$ times.  Frequency permutation arrays were introduced in the 2006 paper of Huczynska and Mullen~\cite{Huczynska06}.  A straightforward example of a frequency permutation array can be obtained by taking a permutation code $\mathcal{C}$ of length $n$ and forming the {\em repetition code}, $\Rep_\lambda(\mathcal{C})$, where each codeword is formed by repeating each codeword of $\mathcal{C}$ $\lambda$ times.  If $\mathcal{C}$ has minimum distance $d$, then clearly $\Rep_\lambda(\mathcal{C})$ has minimum distance $\lambda d$.  

For reasons of improved decoding performance, it is therefore desirable to obtain FPAs with the same length, alphabet and size as $\Rep_\lambda(\mathcal{C})$, but with a larger minimum distance.  An approach to this was introduced in the 2015 paper of Gillespie, Praeger and Spiga~\cite{Gillespie15} and further developed by Akbari, Gillespie and Praeger in 2018~\cite{Akbari18}, where {\em twisted permutation codes} were considered.  Informally, the idea is that instead of repeating the same permutation $\lambda$ times over, a codeword can be formed by taking the image of the same element of an abstract group from multiple permutation representations of the same degree; it transpires that this can result in improved minimum distance.  Formally, these are defined as follows.

\begin{defn} \label{defn:twisted}
Let $G$ be an abstract finite group, and let $\mathcal{I}=(\rho_1,\ldots,\rho_\lambda)$ be a $\lambda$-tuple of (not necessarily distinct) permutation representations of $G$ in the symmetric group $S_n$.  For $g\in G$, let $\rho_i(g)$ be written in list form.  The {\em twisted permutation code}, $\Tw(G,\mathcal{I})$, is defined as
\[ \Tw(G,\mathcal{I}) = \left\{ [\, \rho_1(g) \mid \rho_2(g) \mid \cdots \mid \rho_\lambda(g)\, ] \, : \, g\in G \right\}. \]
\end{defn}
That is, each element of $\Tw(G,\mathcal{I})$ is the concatenation of the images of $g$ under each $\rho_i$ (written in list form), so we have a frequency permutation array with alphabet size $n$ and length $\lambda n$.  We call the subwords $\rho_1(g),\rho_2(g),\ldots,\rho_\lambda(g)$ the {\em components} of a codeword (and similarly, we will refer to the components of a received word).

Unlike~\cite{Akbari18,Gillespie15}, we will insist that each permutation representation in $\mathcal{I}=(\rho_1,\ldots,\rho_\lambda)$ is faithful (although all the examples in~\cite{Akbari18,Gillespie15} are faithful).  In the case where $\rho_1,\ldots,\rho_\lambda$ are all the same permutation representation, then we have the $\lambda$-fold repetition code $\Rep_\lambda(\rho_i(G))$ once again.  For a given $G$ and $\mathcal{I}=(\rho_1,\ldots,\rho_\lambda)$, let $\delta_{\mathrm{rep}}$ be the minimum of all of the minimum distances of the $\lambda$-fold repetition codes $\Rep_\lambda(\rho_i(G))$ (for $1\leq i\leq\lambda$), and let $\delta_{\mathrm{tw}}$ be the minimum distance of $\Tw(G,\mathcal{I})$.  In \cite[Theorem~1.1]{Gillespie15}, it is proved that $\delta_{\mathrm{tw}} \geq \delta_{\mathrm{rep}}$, and a number of examples are given where the inequality is strict, such as the following.

\begin{example} \label{example:outS6}
Let $G$ be the symmetric group $S_6$, and let $\rho_1,\rho_2$ be the distinct permutation representations of $S_6$, interchanged by the outer automorphism.  Now, since $S_6$ has minimum distance~$2$, we have that $\Rep_2(\rho_1(G))$ and $\Rep_2(\rho_2(G))$ both have minimum distance~$4$.  However, $\Tw(G,(\rho_1,\rho_2))$ has minimum distance~$8$ (see~\cite[Section 4.1]{Gillespie15}).
\end{example}

In~\cite{Akbari18,Gillespie15}, a number of examples of twisted permutation codes with improved minimum distance (i.e.\ where $\delta_{\mathrm{tw}} > \delta_{\mathrm{rep}}$) are presented, but no decoding algorithm is given.  The purpose of the present paper is to adapt the approach of~\cite{ecpg} for decoding permutation groups as codes to the newer situation of twisted permutation codes.

\section{General results} \label{sec:general}

The following notion is crucial to the decoding algorithm (for permutation groups) in~\cite{ecpg}, and in what follows.

\begin{defn} \label{defn:base}
Let $G$ be a permutation group acting on a finite set $\Omega$.  A {\em base} for $G$ is a subset $\{x_1,\ldots,x_k\}$ of elements of $\Omega$ whose pointwise stabilizer in $G$ is trivial.  The {\em base size} of $G$, denoted $b(G)$, is the smallest size of a base for $G$.
\end{defn}

A direct consequence of the definition is that the action of an element $g\in G$ on a base uniquely identifies $g$: if $(x_1^g,\ldots,x_k^g)=(x_1^h,\ldots,x_k^h)$ then $g=h$.  These are useful for decoding: if a permutation is transmitted and the received word contains errors outside of the positions labelled by a base, then the transmitted permutation can be identified correctly.  However, as the errors could be in any possible positions, a single base will not be sufficient.  Instead, we have the following definition (also taken from~\cite{ecpg}).

\begin{defn} \label{defn:UBB}
Let $G$ be a permutation group acting on a finite set $\Omega$, and let $r\geq 0$.  An {\em uncovering-by-bases of strength $r$} (or {\em $r$-UBB}) for $G$ is a collection $\mathcal{U}$ of bases for $G$ with the property that any $r$-subset of $\Omega$ is disjoint from at least one base in $\mathcal{U}$.
\end{defn}

If $G$ has minimum distance $d$, then we usually assume that $r=\left\lfloor (d-1)/2 \right\rfloor$, which we call the {\em correction capability} of $G$.  

\begin{example} \label{example:UBB}
Consider the group $G=\mathrm{PGL}(2,7)$ in its action on $\Omega=\{1,\ldots,8\}$.  This action is sharply $3$-transitive, so any $3$-tuple from $\Omega$ forms a base, and the minimum distance is $5$, so the correction capability is $\left\lfloor (5-1)/2 \right\rfloor = 2$.  The following is a $2$-UBB for $G$:
\[ 
\begin{array}{ccc}
1 & 2 & 3 \\
4 & 5 & 6 \\
2 & 3 & 7 \\
1 & 7 & 8
\end{array}
\]
By inspection, we see that any pair from $\{1,\ldots,8\}$ is disjoint from at least one triple in the UBB.
\end{example}

We observe that, if the bases in $\mathcal{U}$ each have size $k$ and each base for is regarded as a $k$-subset, then the complements of the bases in $\mathcal{U}$ form an {\em $(n,n-k,r)$ covering design} (see~\cite[{\S}VI.11]{handbook}).  The online database maintained by Gordon~\cite{Gordon} is a useful resource for examples of covering designs with small parameters.  We also remark that for $r\leq \left\lfloor (d-1)/2 \right\rfloor$, an $r$-UBB is guaranteed to exist (see~\cite[Proposition 7]{ecpg}).  Constructions of UBBs for many permutation groups can be found in~\cite{btubb,ecpg,soc}.

The decoding algorithm given in~\cite{ecpg} works as follows: suppose a permutation $g\in G$ is transmitted and the received word $w$ contains at most $r$ errors.  For each base in $\mathcal{U}$, identify the element (if one exists) of $G$ which agrees with $w$ in the positions labelled by the base; if this permutation is within distance $r$ of $w$ then it must be the transmitted permutation $g$.  Since any combination of $r$ error positions is avoided by at least one base in $\mathcal{U}$, we are guaranteed to succeed.

When we speak of a ``base for a group $G$'', it is a property of the specified permutation representation of $G$.  In general, if $G_1$ and $G_2$ are isomorphic groups acting on the same set $\Omega$, it is not necessarily true that a base for $G_1$ is a base for $G_2$.  However, if the following stronger condition holds, the situation is more straightforward.

\begin{defn} \label{defn:perm-iso}
Let $G_1$ and $G_2$ be groups acting on $\Omega_1$ and $\Omega_2$, respectively, and suppose there is an isomorphism $\varphi\, :\, G_1 \to G_2$.  Then $G_1$ and $G_2$ are {\em permutationally isomorphic} if there is a bijection $\psi\, :\, \Omega_1 \to \Omega_2$ such that $\psi( x^g ) = (\psi(x))^{(g^\varphi)}$ for all $x\in \Omega_1$ and all $g\in G_1$.  The pair $(\psi,\varphi)$ is called a {\em permutational isomorphism}.
\end{defn}

In other words, if $G_1$ and $G_2$ are permutationally isomorphic, then not only are they isomorphic as abstract groups, but they act in the same way on their respective domains $\Omega_1$ and $\Omega_2$.  The next result is a straightforward exercise for the reader.

\begin{prop} \label{prop:base-iso}
Suppose that $G_1$ and $G_2$ are groups acting on $\Omega_1$ and $\Omega_2$, respectively, such that $(\psi,\varphi)$ is a permutational isomorphism.  Then if $B=\{x_1,\ldots,x_b\} \subseteq \Omega_1$ is a base for $G_1$ in its action on $\Omega_1$, then $\psi(B) = \{ \psi(x_1),\ldots,\psi(x_b) \} \subseteq \Omega_2$ is a base for $G_2$ in its action on $\Omega_2$.
\end{prop}

As a consequence, if we have two permutationally-isomorphic groups $G_1$ and $G_2$, we can obtain an uncovering-by-bases for $G_2$ by applying the map $\psi$ to the bases in a UBB for $G_1$.

\subsection{Adapting the algorithm}

To adapt the algorithm from \cite{ecpg} to twisted permutation codes, we recall that the codewords in $\Tw(G,\mathcal{I})$ are in one-to-one correspondence with the elements of the abstract group $G$, so decoding will still involve identifying group elements.  For each $i$, we let $G_i$ denote the image of the faithful representation $\rho_i$, so that $\mathcal{G}=(G_1,\ldots,G_\lambda)$ is a list of permutation groups of degree $n$ isomorphic to $G$.  Without loss of generality, we pick $G_1$ as a ``distinguished'' copy.  For now, assume that each $G_i$ is permutationally isomorphic to $G_1$.

Next, define $\alpha_i\, : \, G_1 \to G_i$ as the composition of $\rho_1^{-1}$ (defined on $G_1={\rm Im}(\rho_1)$) with $\rho_i$.  Consequently, $\mathcal{A}=(\alpha_1,\ldots,\alpha_\lambda)$ gives a list of isomorphisms from $G_1$ to each of $(G_1,\ldots,G_\lambda)$, while $\mathcal{A}^{-1}=(\alpha_1^{-1},\ldots,\alpha_\lambda^{-1})$ gives their respective inverses.  We also let $\mathcal{F}=(\psi_1,\ldots,\psi_\lambda)$ be bijections such that $(\psi_i,\alpha_i)$ is a permutational isomorphism from $G_1$ to $G_i$ (for $1\leq i\leq\lambda$), and $\mathcal{F}^{-1} = (\psi_1^{-1},\ldots,\psi_\lambda^{-1})$ gives the respective inverses.  (Note that in the case of the repetition code $\Rep_\lambda(G_1)$, each $\alpha_i$ and each $\psi_i$ is the identity map.)

Let $r_{\mathrm{tw}}=\left\lfloor (\delta_{\mathrm{tw}}-1)/2 \right\rfloor$ be the correction capability of $\Tw(G,\mathcal{I})$, and let $r'=\left\lfloor r_{\mathrm{tw}}/\lambda \right\rfloor$; by the pigeonhole principle, if a received word contains at most $r_{\mathrm{tw}}$ errors spread across $\lambda$ components, then there must be a component containing at most $r'$ errors.  Finally, suppose that $\mathcal{U}$ is a UBB for $G_1$ of strength $r'$ (so $\psi_i(\mathcal{U}) = \{ \psi_i(B)\, :\, B\in\mathcal{U} \}$ is a UBB for $G_i$ of strength $r'$).

In an implementation of the algorithm, the receiver knows the list of groups $G_1,\ldots,G_\lambda$, as well as the lists of mappings $\mathcal{A}$, $\mathcal{A}^{-1}$, $\mathcal{F}$ and $\mathcal{F}^{-1}$, the correction capability $r_{\mathrm{tw}}$, and the uncovering-by-bases $\mathcal{U}$.  An input to the algorithm consists of a received word $\mathbf{w}=[w_1,\ldots,w_\lambda]$.

\begin{algorithm} \label{algorithm:main}
Suppose that the transmitted codeword is $\mathbf{g}=[\rho_1(g),\ldots,\rho_\lambda(g)]$, and that the received word $\mathbf{w}=[w_1,\ldots,w_\lambda]$ contains at most $r_{\mathrm{tw}}$ errors.  Choose the first base $B_1\in \mathcal{U}$ and examine the symbols in $w_1$ in the positions indexed by $B_1$; if there are no repeated symbols, find (if it exists) the unique element $h_1\in G_1$ which agrees with $w_1$ in those positions, then compute $h_1^{\alpha_2} \in G_2, \ldots, h_1^{\alpha_\lambda} \in G_\lambda$ to obtain a codeword $\mathbf{h}= [h_1,h_1^{\alpha_2},\ldots,h_1^{\alpha_\lambda}] \in \Tw(G,\mathcal{I})$.  If $\hh$ is within distance $r_{\mathrm{tw}}$ of $\ww$, we must have that $\hh=\gg$, and we have decoded successfully.

Otherwise, we move to $w_2$ and examine the symbols in the positions of $w_2$ indexed by the base $\psi_2(B_1)$ for $G_2$; if there are no repeats, find (if it exists) the unique element $h_2\in G_2$ which agrees with $w_2$ in those positions, then compute $h_2^{\alpha_2^{-1}} \in G_1$ as well as $h_2^{\alpha_2^{-1} \alpha_3}\in G_3, \ldots, h_2^{\alpha_2^{-1} \alpha_\lambda}\in G_\lambda$ to obtain a codeword $\hh= [h_2^{\alpha_2^{-1}}, h_2, h_2^{\alpha_2^{-1} \alpha_3}, \ldots, h_2^{\alpha_2^{-1} \alpha_\lambda}] \in \Tw(G,\mathcal{I})$.  Again, if $\hh$ is within distance $r_{\mathrm{tw}}$ of $\ww$, we have decoded successfully.

This process is then continued for each $w_i$ until we can decode successfully; if we fail each time, we then consider the next base $B_2 \in\mathcal{U}$ and repeat the process for each $w_1,\ldots,w_\lambda$, examining the positions in each $w_i$ labelled by $\psi_i(B_2)$, reconstructing a codeword in $\Tw(G,\mathcal{I})$ and comparing it to $\ww$.  If there is still no success, we move to the next base in $\mathcal{U}$, and then the next, and so on, until we are successful.
\end{algorithm}

Since there must be a component $w_i$ which contains at most $r'$ errors, and because $\mathcal{U}$ is an uncovering-by-bases of strength $r'$, there must be a base $\psi_i(B_j)$ for $G_i$ avoiding these errors.  So the algorithm is guaranteed to decode successfully after at most $\lambda|\mathcal{U}|$ attempts.

\begin{figure}[h]
\centering
\begin{tikzpicture}[node distance=1.6cm]

\tikzstyle{startstop} = [rectangle, rounded corners, minimum width=2.4cm, minimum height=0.75cm,text centered, draw=black]
\tikzstyle{process} = [rectangle, minimum width=2.4cm, minimum height=0.75cm, text centered, draw=black]
\tikzstyle{decision} = [diamond, minimum width=2.4cm, minimum height=0.75cm, text centered, draw=black]
\tikzstyle{arrow} = [thick,->,>=stealth]

\node (start) [startstop] {\scriptsize Start};
\node (pro1) [process, below of=start] {\scriptsize Separate received word into components};
\node (pro2) [process, below of=pro1] {\scriptsize Choose first base};
\node (pro3) [process, below of=pro2] {\scriptsize Choose first component};
\node (pro4) [process, below of=pro3] {\scriptsize Identify corresponding group element};
\node (pro5) [process, below of=pro4] {\scriptsize Build codeword using inverse maps};
\node (dec1) [decision, below of=pro5, yshift=-1.6cm,align=left] {\scriptsize Within distance $r_{\mathrm{tw}}$ of \\ \scriptsize received word?};
\node (pro2b) [process, right of=dec1, xshift=2.75cm] {\scriptsize Choose next component};
\node (pro2c) [process, right of=pro2b, xshift=3.5cm] {\scriptsize Choose next base};
\node (stop) [startstop, below of=dec1, yshift=-1.3cm] {\scriptsize Stop};

\draw [arrow] (start) -- (pro1);
\draw [arrow] (pro1) -- (pro2);
\draw [arrow] (pro2) -- (pro3);
\draw [arrow] (pro3) -- (pro4);
\draw [arrow] (pro4) -- (pro5);
\draw [arrow] (pro5) -- (dec1);
\draw [arrow] (dec1) -- node[anchor=east]{\scriptsize Yes}(stop);
\draw [arrow] (dec1) -- node[anchor=south]{\scriptsize No}(pro2b);
\draw [arrow] (pro2b) -- node[anchor=south]{\scriptsize No components}(pro2c);
\draw [arrow] (pro2b) -- node[anchor=north]{\scriptsize remaining}(pro2c);
\draw [arrow] (pro2b) |- (pro4);
\draw [arrow] (pro2c) |- (pro3);

\end{tikzpicture}
\caption{Decoding algorithm for a twisted permutation code}
\label{fig:TwistAlg}
\end{figure}
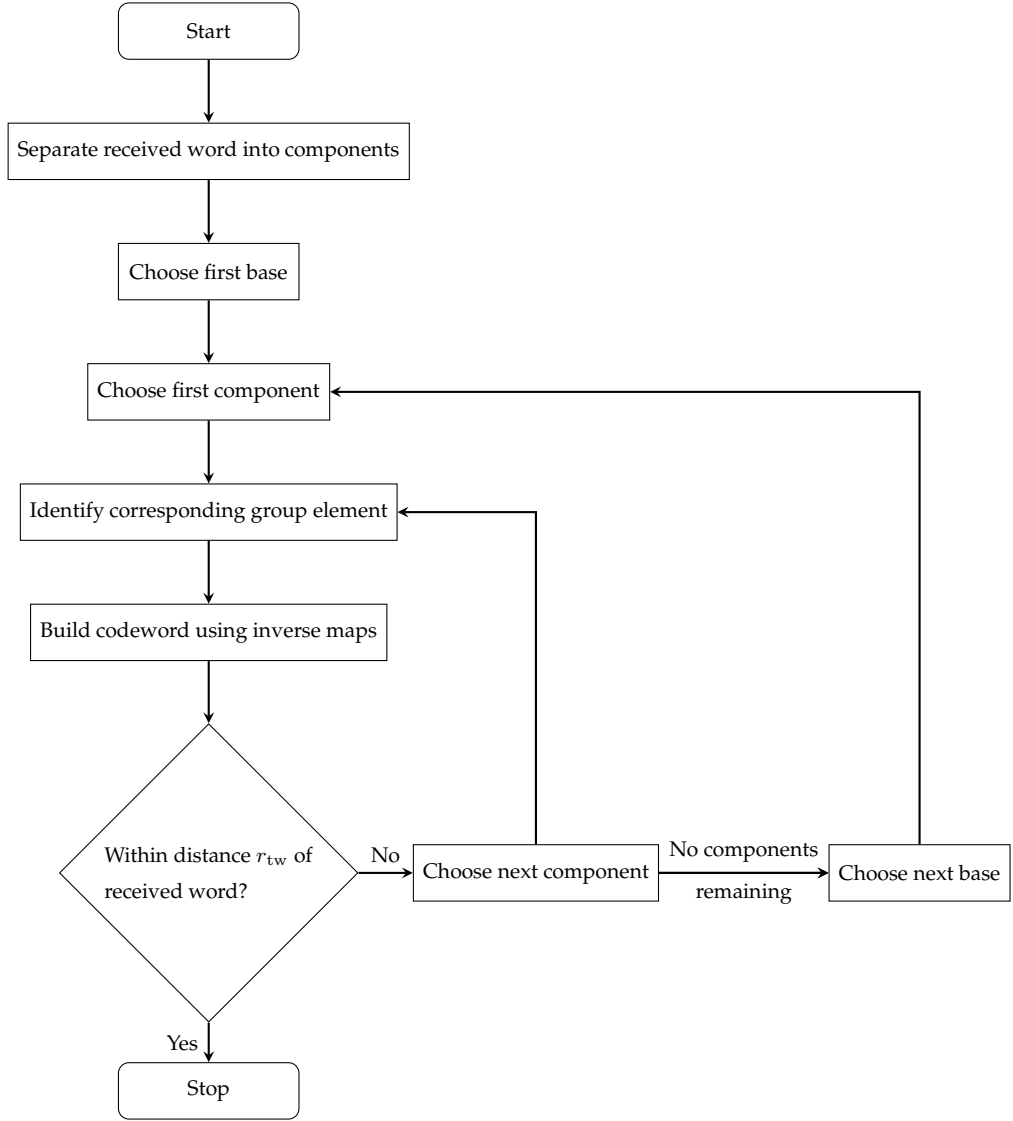

Figure~\ref{fig:TwistAlg} gives a flowchart depicting Algorithm~\ref{algorithm:main}.

\begin{example} \label{example:ASL}
Consider the group $G=\mathrm{ASL}(3,2)$ of affine transformations of $\FF_2^3$.  There are two distinct permutation representations $\rho_1,\rho_2$ as subgroups of $S_8$, interchanged by an outer automorphism; using {\sf GAP}~\cite{GAP}, we find one (in terms of its action on generators of $G$) to be
\begin{eqnarray*}
(2,5)(4,7) & \mapsto & (1,3)(2,7)(4,5)(6,8), \\
(2,3,4)(5,6,8) & \mapsto & (2,3,4)(5,6,8), \\
(1,2)(3,4)(5,6)(7,8) & \mapsto & (1,2)(3,4)(5,6)(7,8), \\
(1,3)(2,4)(5,7)(6,8) & \mapsto & (1,3)(2,4)(5,7)(6,8), \\
(1,5)(2,6)(3,7)(4,8) & \mapsto & (1,5)(2,6)(3,7)(4,8).
\end{eqnarray*}
Denote this mapping by $\alpha_2$, and consider the twisted permutation code $\Tw(G,\mathcal{I})$ where $\mathcal{I}=(\rho_1,\rho_2)$.  We have that $\mathcal{A}=(\id, \alpha_2)$, and since $\im(\rho_1)=\im(\rho_2)$ we have $\mathcal{F}=(\id,\id)$.  Now, by \cite[subsection 7.1]{Gillespie15}, we have $\delta_{\mathrm{tw}}=12$, which is an improvement on $\delta_{\mathrm{rep}}=8$; this means that $\Tw(G,\mathcal{I})$ can correct $\lfloor (12-1)/2\rfloor=5$ errors, while $r'=\lfloor 5/2\rfloor =2$.  

$G$ has base size $4$, with the minimum bases corresponding to affine-independent $4$-tuples in $\FF_2^3$; below is an uncovering-by-bases $\mathcal{U}$ of strength $r'=2$ for $G$:
\begin{center}
\begin{tabular}{ c c c c}
  1 & 2 & 3 & 5\\
  4 & 5 & 6 & 7\\
  1 & 4 & 6 & 8\\
  1 & 5 & 7 & 8\\
  2 & 3 & 4 & 6\\
  2 & 3 & 7 & 8
\end{tabular}
\end{center}
We can use {\sf GAP} to verify that each row is a base for $G$, and by inspection any $2$-subset of $\{1,\ldots,8\}$ is disjoint from at least one base in $\mathcal{U}$.

Let $g=(1,4,6,8,5,3)(2,7)\in G$.  Now, $g^{\alpha_2} = (1,7,5,8,2,4)(3,6)$, so by concatenating these in list form we obtain the codeword
\[ \gg = [ g_1 \mid g_2 ] = [ 4, 7, 1, 6, 3, 8, 2, 5 \mid 7, 4, 6, 1, 8, 3, 5, 2 ] \] 
in $\Tw(G,\mathcal{I})$.  Suppose that $\gg$ is transmitted, and the following word (with two errors) is received:
\[ \ww = [ w_1 \mid w_2 ] = [ 4, 7, 1, 6, 7, 8, 2, 5 \mid 4, 4, 6, 1, 8, 3, 5, 2 ]. \]

The first base is $\{1,2,3,5\}$, so we first examine the symbols in those positions of $w_1$, which are $4,7,1,7$; since symbol $7$ is repeated, we cannot decode. We then examine those positions of $w_2$ and find $4,4,6,8$, so we are unsuccessful again.

The next base is $\{4,5,6,7\}$; in $w_1$ we find symbols $6,7,8,2$. There are no repeats, so we obtain the element $h_1=[4,3,5,6,7,8,2,1]\in G_1$ which agrees with $w_1$ in those positions; applying $\alpha_2$ yields $h_1^{\alpha_2}= [1, 2, 8, 7, 6, 5, 3, 4]$, and we obtain the codeword
\[ \hh = [ 4, 3, 5, 6, 7, 8, 2, 1 \mid 1, 2, 8, 7, 6, 5, 3, 4 ]. \]
However, since this is at distance $11>\delta_{\mathrm{tw}}$ from $\ww$, it is rejected.  In $w_2$, we find symbols $1,8,3,5$, and obtain $h_2=[ 7, 4, 6, 1, 8, 3, 5, 2 ]$; applying $\alpha_2^{-1}$ yields $h_2^{\alpha_2^{-1}} = [ 4, 7, 1, 6, 3, 8, 2, 5 ]$, and we obtain the codeword
\[ \hh = [ 4, 7, 1, 6, 3, 8, 2, 5 \mid 7, 4, 6, 1, 8, 3, 5, 2 ] \]
which is distance $2$ from $\ww$.  So we can conclude that $\hh=\gg$, and we have decoded successfully.
\end{example}

In the more general case where the groups $G_1,\ldots,G_\lambda$ are not permutationally isomorphic, we no longer have the list of mappings $\psi_1,\ldots,\psi_\lambda$, and will require a separate uncovering-by-bases  for each distinct image group $G_i$.  The decoding algorithm proceeds similarly, but the need for additional UBBs makes it more difficult to implement.  However, in each of the examples considered in \cite{Akbari18,Gillespie15}, the image groups $G_i$ are typically not just permutationally isomorphic, but are in fact equal.  This means that each map $\psi_i$ is the identity map, and we can use the same UBB in each component.

\subsection{Decoding repetition codes and ``unimproved'' codes} \label{sec:repetition}

The following observation is helpful, as it ensures that known UBBs for permutation codes can be applied to repetition codes and ``unimproved'' twisted permutation codes (i.e.\ those for which $\delta_{\mathrm{tw}}=\delta_{\mathrm{rep}}$).

\begin{prop} \label{prop:repetition-strength}
Suppose that $G$ is a permutation group with correction capability $r$.  Then, for the repetition code $\Rep_\lambda(G)$, the strength $r'$ of the UBB required for $\Rep_\lambda(G)$ is equal to $r$.
\end{prop}

\proof We know that $r=\lfloor (d-1)/2\rfloor$, where $d$ is the minimum distance of $G$.  Now, $\Rep_\lambda(G)$ has minimum distance $\lambda d$, correction capability $\lfloor (\lambda d-1)/2 \rfloor$, and we have $r'=\lfloor \lfloor (\lambda d-1)/2 \rfloor /\lambda \rfloor$.  A case analysis to consider when $d$ and $\lambda$ are each odd or even then shows that, in all cases, $r'=\lfloor (\lambda d-1)/2 \rfloor=r$. \endproof

Note that the same result holds for ``unimproved'' twisted permutation codes; consequently, the UBBs obtained in \cite{btubb,ecpg,soc} may be used not just for the corresponding repetition codes, but also the ``unimproved'' twisted permutation codes.  In the table below, we give some further examples of groups (with multiple permutation representations), including some mentioned in~\cite{Gillespie15} (namely $\mathrm{PSL}(2,11)$, $M_{12}$ and $A_7$), and their parameters.  For $2^4:A_6$, $2^4:S_6$ and $M_{22}$, we verified that $\delta_{\mathrm{tw}} = \delta_{\mathrm{rep}}$ with the same techniques as~\cite{Gillespie15}, using {\sf GAP}.

\begin{table}[hbtp]
\centering
\renewcommand{\arraystretch}{1.1}
\begin{tabular}{c|cccccccc} \hline
Group $G$   & $|G|$    & $n$ & $\lambda$ & $\delta_{\mathrm{tw}} = \delta_{\mathrm{rep}}$ 
                                                            & $r$  & $r'$ & $b(G)$ & $|\mathcal{U}|$ \\ \hline
$\mathrm{PSL}(2,11)$ & $660$    & $11$ & $2$      & $2\cdot 8 = 16$  & $7$  & $3$  & $3$    & $5$ \\
$M_{12}$             & $95040$  & $12$ & $2$      & $2\cdot 8 = 16$  & $7$  & $3$  & $5$    & $11$ \\
$A_7$                & $2520$   & $15$ & $2$      & $2\cdot 12 = 24$ & $11$ & $5$  & $3$    & $9$ \\
$2^4:A_6$            & $5760$   & $16$ & $4$      & $4\cdot 12 = 48$ & $23$ & $5$  & $4$    & $12$ \\
$2^4:S_6$            & $11520$  & $16$ & $2$      & $2\cdot 8 = 16$  & $7$  & $3$  & $5$    & $6$ \\
$M_{22}$             & $443520$ & $22$ & $2$      & $2\cdot 16 = 32$ & $15$ & $7$  & $5$    & $22$ \\ \hline
\end{tabular}
\renewcommand{\arraystretch}{1.0}
\caption{Parameters and decoding for some ``unimproved'' twisted permutation codes}
\label{table:unimproved}
\end{table}
The UBBs mentioned in Table~\ref{table:unimproved} can be found in Appendix~\ref{app:UBBs}.  Each was obtained by taking the complements of blocks of the corresponding $(n,n-b(G),r')$-covering designs given in Gordon's database~\cite{Gordon}; in some cases, the points needed to be relabelled to ensure that the complement of each block was a base for the group $G$.

\subsection{Decoding twisted permutation codes with improved minimum distance} \label{sec:improved-UBBs}

In the case of a twisted permutation code $\Tw(G,\mathcal{I})$ with improved minimum distance, the strength of the UBB we need is typically larger than that needed for the repetition code $\Rep_\lambda(G)$.  For the groups $S_6$, $A_6$ and $\mathrm{ASL}(3,2)$, each of which were shown in~\cite{Gillespie15} to yield such ``improved'' codes, we summarize the details in Table~\ref{table:improved} below.

\begin{table}[hbtp]
\centering
\renewcommand{\arraystretch}{1.1}
\begin{tabular}{c|ccccccccc} \hline
Group $G$            & $|G|$    & $n$ & $\lambda$ & $\delta_{\mathrm{rep}}$ & $\delta_{\mathrm{tw}}$
                                                      & $r_{\mathrm{tw}}$ & $r'$ & $b(G)$ & $|\mathcal{U}|$ \\ \hline
$S_6$                & $720$    & $6$ & $2$       & $4$ & $8$  & $3$               & $1$  & $5$    & $6$ \\
$A_6$                & $360$    & $6$ & $2$       & $6$ & $8$  & $3$               & $1$  & $4$    & $3$ \\
$\mathrm{ASL}(3,2)$  & $1344$   & $8$ & $2$       & $8$ & $12$ & $5$               & $2$  & $4$    & $6$ \\ \hline
\end{tabular}
\renewcommand{\arraystretch}{1.0}
\caption{Parameters and decoding for some ``improved'' twisted permutation codes}
\label{table:improved}
\end{table}
For $S_6$, the corresponding UBB consists of all $5$-subsets of $\{1,\ldots,6\}$; for $A_6$, we can use $\{ 1234, \linebreak 1256, 3456 \}$; for $\mathrm{ASL}(3,2)$, the UBB is given in Example~\ref{example:ASL}.

As Table~\ref{table:improved} does not give very many examples, in the next section we consider an infinite family of ``improved'' codes, taken from~\cite{Akbari18}.


\section{Codes from the groups $G_k(p)$} \label{sec:Gkp}

In~\cite[Section 3]{Akbari18}, Akbari {\em et al.}\ give an infinite family of twisted permutation codes, which arise from affine groups over the vector space $\FF_p^k$.  We begin by giving a summary of these groups and the corresponding codes.

Let $B_k$ be the following $k\times k$ matrix over $\FF_p$:
\[ B_k = \begin{bmatrix}
         1 & 0 & 0 & \cdots & 0 \\
				 1 & 1 & 0 & \cdots & 0 \\
				 0 & 1 & 1 & \cdots & 0 \\
				 \vdots && \ddots & \ddots & \vdots \\
				 0 & \cdots & 0 & 1 & 1
				\end{bmatrix}. \]
Since $B_k$ is lower unitriangular, it is clearly invertible.  It can be shown (see~\cite[Lemma 4]{Akbari18}) that $B_k$ has multiplicative order $p$, and that its powers are given by
\[ B_k^i = \begin{bmatrix}
				 1 & 0 & 0 & \cdots & 0 \\
				 i & 1 & 0 & \cdots & 0 \\
				 \binom{i}{2} & i & 1 & \cdots & 0 \\
				 \vdots & \vdots& \ddots & \ddots & \vdots \\
				 \binom{i}{k-1} & \binom{i}{k-2} & \binom{i}{k-3}  & \cdots & 1
				\end{bmatrix} \]
for $0\leq i < p$, and where the entries $\binom{i}{j}$ are taken modulo~$p$.

\begin{defn} \label{defn:Gkp}
Suppose that $V=\FF_p^k$ is the space of row vectors, and let $H_k(p)$ be the subgroup of $GL(k,p)$ generated by $B_k$.  Denote by $G_k(p)$ the subgroup $V\rtimes H_k(p)$ of $\mathrm{AGL}(k,p)$.
\end{defn}

Since $B_k$ has multiplicative order $p$, it follows that $G_k(p)$ has order $p^{k+1}$.  Next, we define a \linebreak $(k+1)\times (k+1)$ matrix
\[ A_{\vv,i} = \begin{bmatrix} 1 & \vv \\
                               \mathbf{0}^T & B_k^i
															 \end{bmatrix}, \]
where $\vv$ is a row vector in $\FF_p^k$, and $0\leq i < p$.  Then denote the collection of all such matrices by $\overline{G_k(p)}$.  In~\cite[Lemma 5]{Akbari18}, it is shown that $\overline{G_k(p)}$ is a subgroup of $GL(k+1,p)$ under matrix multiplication, has order $p^{k+1}$, and is isomorphic to the affine group $G_k(p)$.  They also show (in~\cite[Lemma 7]{Akbari18}) that $\overline{G_k(p)}$ acts faithfully and transitively on the set of $p^k$ row vectors $\Omega =\{ (1,\vv) \,:\, \vv\in\FF_k^p \}$, so $\overline{G_k(p)}$ can be viewed as a permutation group on $p^k$ points.  Furthermore, they give a collection of isomorphisms, which they denote by $\tau_\ww$, from $G_k(p)$ to $\overline{G_k(p)}$; these give a collection of permutation representations of $G_k(p)$ in $\Sym(\Omega)$, and thus can be used to construct twisted permutation codes.  In particular, they show that as a permutation code $\overline{G_k(p)}$ has minimum distance $p^k-p$, so the $p$-fold repetition code $\Rep_p(\overline{G_k(p)})$ has minimum distance $p^{k+1}-p^2$; however, using a particular collection $\mathcal{I}$ of $p$ permutation representations, there is a twisted permutation code $\Tw(G_k(p), \mathcal{I})$ with improved minimum distance $p^{k+1}-p$ (see~\cite[Proposition 9]{Akbari18}).

To apply Algorithm~\ref{algorithm:main} to these codes, we first need to obtain bases for the group $\overline{G_k(p)}$.  Let $\ee_1,\ldots,\ee_k$ denote the standard basis vectors of $\FF_p^k$.

\begin{prop} \label{prop:gkp-base}
For any $j$ where $2\leq j\leq k$, we have that $\{ (1,\mathbf{0}),\, (1,\ee_j) \}$ is a base for $\overline{G_k(p)}$ acting on the set $\Omega$.  Furthermore, $b(\overline{G_k(p)})=2$.
\end{prop}

\proof Suppose that $A_{\vv,i}$ lies in the pointwise stabilizer of $\{ (1,\mathbf{0}),\, (1,\ee_j) \}$ in $\overline{G_k(p)}$.  First, we have
\begin{eqnarray*}
(1,\mathbf{0})A_{\vv,i} & = & (1,\mathbf{0})\begin{bmatrix} 1 & \vv \\ \mathbf{0}^T & B_k^i \end{bmatrix}  \\
                        & = & (1+0, \vv+\mathbf{0}) \\
												& = & (1,\vv),
\end{eqnarray*}
so for $A_{\vv,i}$ to fix $(1,\mathbf{0})$ we must have $\vv=\mathbf{0}$.  Then we have
\begin{eqnarray*}
(1,\ee_j)A_{\mathbf{0},i} & = & (1,\ee_j)\begin{bmatrix} 1 & \mathbf{0} \\ \mathbf{0}^T & B_k^i \end{bmatrix} \\
                          & = & (1+0, \mathbf{0}+\ee_j B_k^i ).
\end{eqnarray*}
Now, $\ee_j B_k^i$ is precisely row $j$ of $B_k^i$, i.e.
\[ \ee_j B_k^i = \left( \binom{i}{j-1}, \binom{i}{j-2}, \cdots, \binom{i}{1}, 1, 0, \cdots 0 \right). \] 
So for this to be equal to $\ee_j$, we require that $j\geq 2$ and that all the binomial coefficients are equal to $0$; this will happen only when $i=0$.  Consequently, we have that $A_{\mathbf{0},i} = A_{\mathbf{0},0} = I_{k+1}$, i.e.\ the identity element of $\overline{G_k(p)}$.  Hence $\{ (1,\mathbf{0}),\, (1,\ee_j) \}$ is a base for $\overline{G_k(p)}$ acting on $\Omega$.

Since $\overline{G_k(p)}$ acts transitively on $\Omega$ but $|\Omega|>|\overline{G_k(p)}|$, it follows that $\overline{G_k(p)}$ has no base of size~$1$; therefore $b(\overline{G_k(p)})=2$.
\endproof

The next step is to obtain an uncovering-by-bases for $\overline{G_k(p)}$.  Now, since $b(\overline{G_k(p)})=2$, we can regard the minimum bases as the edges of a graph with vertex set $\Omega$.  The following terminology and notation was introduced by in 2020 by Burness and Giudici~\cite{Burness20}.

\begin{defn} \label{defn:saxlgraph}
Let $G$ be a group acting on $\Omega$ with $b(G)=2$.  The {\em Saxl graph} of $G$, denoted $\Sigma(G)$, is the graph with vertex set $\Omega$, and where $\{u,v\}$ is an edge if and only if it is a base for $G$.
\end{defn}

Similar graphs (named ``base-orbital graphs'') appear in~\cite[Section 3]{soc}, but where the edge set consists of a single orbit of $G$ on its bases of size~$2$, rather than all such bases.  Thus the edge set of the Saxl graph $\Sigma(G)$ is the union of the edge sets of each of the base-orbital graphs of $G$.

\begin{prop} \label{prop:saxl-matching}
Let $G$ be a permutation group, acting on a set $\Omega$ of size $n$, with base size $b(G)=2$.  Then an optimal uncovering-by-bases for $G$ is a matching in the Saxl graph $\Sigma(G)$.
\end{prop}

\proof The smallest possible size of a UBB of strength $r'$ is $r'+1$, as otherwise there will be a set of $r'$ points which intersects each base non-trivially. A set of $r'+1$ disjoint bases (i.e.\ a set of $r'+1$ disjoint edges in $\Sigma(G)$) will be sufficient.  But this is exactly a matching of size $r'+1$ in $\Sigma(G)$.  Since $r' \leq \left\lfloor \frac{n-2}{2} \right\rfloor$, the requirement that $2(r+1)\leq n$ will always hold.
\endproof

Recall that a {\em perfect matching} in a graph $\Gamma$ on $n$ vertices is a matching (of size $\frac{n}{2}$) which includes every vertex of $\Gamma$.  If $n$ is odd, no perfect matching can exist, but a {\em near-perfect matching} is a matching which includes every vertex except one.  The following result is well-known (see, for example, Godsil and Royle~\cite[Section 3.5]{Godsil01}).

\begin{lemma} \label{lemma:matching}
Let $\Gamma$ be a connected, vertex-transitive graph with $n$ vertices. Then $\Gamma$ has either a perfect matching or a near-perfect matching, depending on the parity of $n$.
\end{lemma}

By construction, the Saxl graph $\Sigma(G)$ will be vertex-transitive whenever $G$ is transitive; in order to apply Lemma~\ref{lemma:matching} to $\Sigma(G)$, one must show that it is connected.  (The condition is necessary: a disconnected graph where all components have odd size can never have a (near-) perfect matching.)  We would like to apply this to the Saxl graph of $\overline{G_k(p)}$.

\begin{lemma} \label{lemma:saxl-connected}
Let $G$ be the group $\overline{G_k(p)}$ acting on the set $\Omega = \{ (1,\vv) \,:\, \vv\in\FF_p^k \}$. Then the Saxl graph $\Sigma(G)$ is connected.
\end{lemma}

\proof We will show that for each element $(1,\vv)\in\Omega$, there exists a path in $\Sigma(G)$ to $(1,\mathbf{0})$.

We saw in Proposition~\ref{prop:gkp-base} that $\{ (1,\mathbf{0}), (1,\ee_j) \}$ is a base for $G$ for $2\leq j\leq k$.  Now consider the orbit of $G$ on such bases; for $A_{\vv,i}\in G$, we have that
\[ (1,\mathbf{0})A_{\vv,i} = (1,\vv) \]
and
\[ (1,\ee_j)A_{\vv,i} = (1, \vv+\bb_j^i) \]
where $\bb_j^i = \ee_j B_k^i$ denotes row $j$ of the matrix $B_k^i$.  So $\{ (1,\vv), (1, \vv+\bb_j^i) \}$ is a base for $G$.  Therefore, in $\Sigma(G)$ each vertex $(1,\vv)$ is adjacent to $(1,\vv\pm \bb_j^i)$; we can label these edges by $\bb_j^i$.  Now, for $2\leq j\leq k$, we have that $\bb_j^i=\ee_j$, while $\bb_2^1 = \ee_1+\ee_2$.

Now suppose that $\vv=(v_1,v_2,\ldots,v_k)$.  Then there is a path in $\Sigma(G)$ from $(1,v_1,v_2,\ldots,v_k)$ to $(1,v_1,v_2,\ldots,v_{k-1},0)$ using edges labelled by $\ee_k$, then a path from $(1,v_1,v_2,\ldots,v_{k-1},0)$ to \linebreak $(1,v_1,v_2,\ldots,v_{k-2},0,0)$ using edges labelled by $\ee_{k-1}$, and so on, until we reach $(1,v_1,v_2,0,\ldots,0)$.  From there, there is a path to $(1,v_1,v_1,0,\ldots,0)$ using edges labelled by $\ee_2$, and then finally a path to $(1,0,\ldots,0)$ using edges labelled by $\ee_1+\ee_2$.

Since there is a path in $\Sigma(G)$ from any vertex to $(1,\mathbf{0})$, it follows that $\Sigma(G)$ is connected.
\endproof

Putting all of this together, we have the following result.

\begin{thm} \label{thm:gkp}
The twisted permutation code $\Tw(G_k(p),\mathcal{I})$, which has size $p^{k+2}$, length $p^{k+1}$ and minimum distance $p^{k+1}-p$, can be decoded using an uncovering-by bases of optimal size $r'+1 = \left\lfloor \frac{p^k-1}{2} \right\rfloor$.
\end{thm}

\proof The size, length and minimum distance of $\Tw(G_k(p),\mathcal{I})$ were all determined in~\cite{Akbari18}.  Calculating $r'$ is a straightforward exercise from this. Since $\Sigma(\overline{G_k(p)})$ is connected (by Lemma~\ref{lemma:saxl-connected}), it has a (near-) perfect matching (by Lemma~\ref{lemma:matching}), which forms an optimal UBB (by Proposition~\ref{prop:saxl-matching}).  \endproof

\section{Another infinite family} \label{sec:ASL2}
We conclude the paper by mentioning another family of groups, mentioned in~\cite[Section~6]{Gillespie15}, with multiple permutation representations.  Let $V$ be the additive group of $\FF_{2^m}^2$ and $K$ denote the special linear group $\mathrm{SL}(2,2^m)$. Since the first cohomology group $H^1(K,V)$ has order $2^m$ (cf.~\cite[Table 7.3]{Cameron99}), we have $2^m$ outer automorphisms of the semidirect product $G=V\rtimes K$.  It follows that there are $2^m$ permutation representations of $G$, the affine special linear group $\mathrm{ASL}(2,2^m)$, and thus $\mathrm{ASL}(2,2^m)$ is a candidate for constructing a twisted permutation code.  Furthermore, a family of UBBs which can be used for these groups is constructed in~\cite[Theorem 5.27]{btubb}.

Unfortunately, as shown in~\cite[Theorem~6.1]{Gillespie15}, twisting does not yield codes with improved minimum distance in this instance.  The possible numbers of fixed points of a non-identity element of $G=\mathrm{ASL}(2,2^m)$ are $0$, $1$ or $2^m$, and thus the minimum distance is $2^{2m}-2^m$.  For the minimum distance of a twisted permutation code $\Tw(G,\mathcal{I})$ to be improved from that of the repetition code $\Rep_\lambda(G)$, a necessary condition is that each conjugacy class of elements of $G$ with the maximum number of fixed points must be mapped by an outer automorphism to a conjugacy class of elements with fewer fixed points; these two conjugacy classes must have the same size and consist of elements of the same order.  However, in $\mathrm{ASL}(2,2^m)$ no such conjugacy classes can exist: there is a unique conjugacy class of elements with $2^m$ fixed points, and these elements have order~$2$; there is only one other conjugacy class of elements of order~$2$, formed of the non-zero elements of the elementary abelian subgroup $V$, and this has a different size.


\appendix

\section{Appendix: Some examples of uncoverings-by-bases} \label{app:UBBs}

\begin{example} \label{example:psl211}
For the group $G=\mathrm{PSL}(2,11)$, we have $n=11$, $b(G)=3$ and $r'=3$:
\begin{center}
\begin{tabular}{c c c}
1 & 2 & 11\\
2 & 8 & 10\\
3 & 4 & 5\\
6 & 7 & 9\\
8 & 10 & 11
\end{tabular}
\end{center}
\end{example}


\begin{example} \label{example:m12}
For the group $M_{12}$, we have $n=12$, $b(M_{12})=5$ and $r'=3$. Since $M_{12}$ is sharply $5$-transitive, any $5$-subset forms a base.  This example is taken from~\cite[Table 1]{btubb}.

\begin{center}
\begin{tabular}{c c c c c}
1 & 2 & 3 & 4 & 5 \\
1 & 2 & 6 & 11 & 12 \\
1 & 3 & 7 & 8 & 9 \\
1 & 4 & 6 & 7 & 10 \\
1 & 5 & 8 & 9 & 11 \\
2 & 4 & 8 & 9 & 12 \\
2 & 5 & 7 & 10 & 11 \\
3 & 4 & 7 & 11 & 12 \\
3 & 5 & 6 & 10 & 12 \\
3 & 6 & 8 & 9 & 11 \\
6 & 7 & 8 & 9 & 10  
\end{tabular}
\end{center}
\end{example}


\begin{example} \label{example:a7}
For the group $A_7$, we have $n=15$, $b(A_7)=4$ and $r'=5$. This example is referred to in~\cite{btubb}, but is not given there explicitly.
\begin{center}
\begin{tabular}{c c c}
1 & 2 & 8\\
2 & 6 & 7\\
3 & 4 & 5\\
6 & 7 & 8\\
9 & 12 & 15\\
9 & 13 & 14\\
10 & 11 & 12\\
10 & 11 & 15\\
13 & 14 & 15
\end{tabular}
\end{center}
\end{example}

\newpage

\begin{example} \label{example:2_4_a6}
For the group $G=2^4:A_6$, we have $n=16$, $b(G)=4$ and $r'=5$:
\begin{center}
\begin{tabular}{c c c c}
1 & 2 & 9 & 16\\
1 & 3 & 9 & 10\\
1 & 4 & 9 & 11\\
2 & 3 & 10 & 16\\
2 & 4 & 11 & 16\\
3 & 4 & 10 & 11\\
5 & 6 & 12 & 13\\
5 & 7 & 13 & 14\\
5 & 8 & 13 &15\\
6 & 7 & 12 & 14\\
6 & 8 & 12 & 15\\
7 & 8 & 14 & 15
\end{tabular}
\end{center}
\end{example}

\begin{example} \label{example:2_4_s6}
For the group $G=2^4:S_6$, we have $n=16$, $b(G)=5$ and $r'=3$:
\begin{center}
\begin{tabular}{c c c c c}
1 & 3 & 4 & 6 & 13 \\ 
2 & 3 & 6 & 11 & 16 \\
2 & 4 & 11 & 13 & 16 \\
5 & 7 & 9 & 14 & 15 \\
7 & 8 & 10 & 12 & 14 \\
8 & 9 & 10 & 12 & 15
\end{tabular}
\end{center}
\end{example}

\begin{example} \label{example:m22}
For the group $M_{22}$, we have $n=22$, $b(M_{22})=5$ and $r'=7$:
\begin{center}
\begin{tabular}{c c c c c}
1 & 2 & 4 & 14 & 19 \\
1 & 2 & 7 & 13 & 17 \\
1 & 3 & 4 & 15 & 17 \\
1 & 3 & 5 & 7 & 19 \\
1 & 5 & 13 & 14 & 15 \\
2 & 3 & 5 & 14 & 17 \\
2 & 3 & 13 & 15 & 19 \\
2 & 4 & 5 & 7 & 15 \\
3 & 4 & 7 & 13 & 14 \\
4 & 5 & 13 & 17 & 19 \\
6 & 8 & 9 & 11 & 18 \\
6 & 8 & 10 & 20 & 21 \\
6 & 9 & 12 & 16 & 21 \\
6 & 10 & 12 & 18 & 22 \\
6 & 11 & 16 & 20 & 22 \\
7 & 14 & 15 & 17 & 19 \\
8 & 9 & 10 & 16 & 22 \\
8 & 11 & 12 & 21 & 22 \\
8 & 12 & 16 & 18 & 20 \\
9 & 10 & 11 & 12 & 20 \\
9 & 18 & 20 & 21 & 22 \\
10 & 11 & 16 & 18 & 21
\end{tabular}
\end{center}
\end{example}

\section*{Acknowledgements}
The first author acknowledges financial support from an NSERC Discovery Grant.  The second author was the recipient of an NSERC Undergraduate Student Research Award.  Part of this work is taken from the second author's B.Sc.\ Senior Project thesis at Grenfell Campus, Memorial University.  Both authors would like to thank the anonymous referees for their suggestions, and Peter Cameron and Robert Guralnick for useful discussions.


\end{document}